\documentclass[aps,pre,twocolumn,groupedaddress]{revtex4-1}

\usepackage{graphicx,color,mathtools,amsmath,amssymb}

\graphicspath{ {./} }

\begin{document}


\title{Data-driven discovery of partial differential equation models with latent variables}


\author{Patrick A.K. Reinbold}
\author{Roman O. Grigoriev}
\affiliation{School of Physics, Georgia Institute of Technology, Atlanta, Georgia 30332-0430, USA}

\date{\today}

\begin{abstract}
In spatially extended systems, it is common to find latent variables that are hard, or even impossible, to measure
with acceptable precision, but are crucially important for the proper description of the dynamics. 
This substantially complicates construction of an accurate model for such systems using data-driven approaches. 
The present paper illustrates how physical constraints can be employed to overcome this limitation using the example of a weakly turbulent quasi-two-dimensional Kolmogorov flow driven by a steady Lorenz force with an unknown spatial profile. Specifically, the terms involving latent variables in the partial differential equations governing the dynamics can be eliminated at the expense of raising the order of that equation.
We show that local polynomial interpolation combined with symbolic regression can handle sparse data on grids that are representative of typical experimental measurement techniques such as particle image velocimetry. 
However, we also find that the reconstructed model is sensitive to measurement noise and trace this sensitivity to the presence of high order spatial and/or temporal derivatives.
\end{abstract}


\maketitle

\section{Introduction}
 
Due to advances in data acquisition, storage, and computational power, data-driven discovery of mathematical models of physical systems, which relies primarily on the empirical observations, has emerged as a viable alternative to more traditional approaches based on, e.g., first-principles derivation. 
While methods for constructing linear models of dynamical systems are very well established 
\cite{aastrom1971}, the progress in model discovery for nonlinear processes is relatively recent, with the earliest efforts focusing on nonlinear ordinary differential or difference equation models of low-dimensional dynamics \cite{crutchfield1987,bongard2007,chou2009,brunton2016}.
Some progress has also been made in model discovery for spatially-distributed systems described by nonlinear partial differential equations (PDEs) \cite{rudy2017}, where all of the state variables are directly observable.

In many cases of practical interest, however, one might be interested in determining the model equations using only some of the state variables.
For instance, the primitive variable description of fluid flows relies on two physical fields: velocity and pressure, only the first of which can typically be measured in experiment with meaningful accuracy.
The presence of latent variables such as pressure makes data-driven model discovery substantially more complicated or incomplete, since existing approaches such as sparse regression \cite{brunton2016}, crucially rely on direct measurements of every state variable which appears in the model.
In particular, using velocity measurements alone only allows reconstruction of the vorticity equation \cite{rudy2017} which describes the evolution of the curl of the velocity, but not its individual components.  

Another common problem in data-driven model discovery is sensitivity to noise in the data.
It is especially acute for spatially distributed systems due to the difficulty in accurately estimating spatial derivatives using sparse noisy data.
As an example, adding just 1\% noise to the data causes errors in the model parameters of order 10\% for the nonlinear Schr\'dinger, KdV, and the vorticity equations, 50\% for the Kuramoto-Sivashinsky equation, and introduced spurious terms in the $\lambda-\omega$ model \cite{rudy2017}.
It is not currently understood on a quantitative level what the impact of noise 
is on the accuracy of the model reconstruction, however.
Neither is it clear how the accuracy of the model reconstruction based on sparse, noisy data can be quantified in the absence of some sort of a reference.

This work uses a representative example of a fluid flow to address several of the open questions, mainly (1) how can we get around the lack of direct measurements of latent variables and (2) how can we quantify the accuracy of the resulting model when the measurements of observable variables are sparse and noisy?
The structure of the paper is as follows. Section \ref{sec:ps} describes the physical problem the model of which we are trying to construct 
and the relevant physical constraints. Section \ref{sec:sr} describes our symbolic regression approach. Section \ref{sec:pi} discusses the use of polynomial approximations for estimating spatial and temporal derivatives. The results are presented in Sect. \ref{sec:r} and our conclusions in Sect. \ref{sec:c}.

\section{Problem statement}
\label{sec:ps}

We will focus on the Lorenz-force-driven flow in a thin electrolyte layer supported by a stationary bottom plate \cite{suri_2014,tithof_2017} as it provides an excellent illustration of the challenges in data-driven discovery of a model for a spatially distributed system with latent variables.
The basic physics and symmetry of the problem imposes a number of constraints on the form of the model and the choice of the fundamental variables.
Being a fluid flow, it is described by two fields, velocity ${\bf u}$ and pressure $p$, so we would expect the dynamics of the fluid flow to be described by evolution equations of the general form
\begin{align}\label{eq:up}
\partial_t{\bf u} &= {\bf N}_u({\bf u},p),\nonumber\\
\partial_tp &= N_p({\bf u},p),
\end{align}
where ${\bf N}_u$ and $N_p$ are some (generally nonlinear) differential operators.

The form of these operators can be constrained by both the physics and the symmetry of the problem; we will start with the latter.
In particular, in order to preserve the rotational symmetry, ${\bf N}_u$ has to be a vector, and so it can be constructed as a linear superposition of terms each of which is a vector. 
Since the fluid layer is thin, the vertical component of the velocity is small compared to the horizontal component and we can consider ${\bf u}$ to be two-dimensional (we can think of ${\bf u}$ as describing the flow at the free surface of the electrolyte).
Furthermore, again due to the small thickness of the fluid layer, both ${\bf u}$ and $p$ can be considered functions of horizontal coordinates $x$ and $y$ and time $t$, but not the vertical coordinate $z$.

There are several ways to construct a vector out of ${\bf u}$, $p$, the gradient operator $\nabla$, and the external forcing field ${\bf f}$ (assumed to be time-independent).
The gravitational acceleration ${\bf g}$, the only other vector quantity in the problem, cannot be included in the two-dimensional model, since the latter does not explicitly include the vertical direction.  
Using one vector object, we can construct three vector fields that are linear in ${\bf u}$, $p$, and ${\bf f}$: $\nabla p$, ${\bf u}$, and ${\bf f}$.
More complicated vector fields can be constructed using powers of $\nabla$ and/or nonlinear functions of $p$, ${\bf u}$, and ${\bf f}$.
We will only consider terms that are linear in $p$ and ${\bf f}$, since $\nabla p$ and ${\bf f}$ both describe the (volumetric) force density and they are linearly related to the time rate of change of the momentum density $\rho\partial_t{\bf u}$ according to Newton's 2nd law.

More vector fields can be constructed using several copies of ${\bf u}$ and $\nabla$.
Keeping terms up to third order in ${\bf u}$ and second order in $\nabla$, we obtain the following evolution equation for the velocity field
\begin{align}\label{eq:u}
\partial_t{\bf u} &= c_1({\bf u}\cdot\nabla){\bf u}+c_2\nabla^2{\bf u}
+c_3{\bf u}+c_4(\nabla\cdot{\bf u}){\bf u}\nonumber\\
&+c_5(\nabla\cdot{\bf u})^2{\bf u}
+c_6(\nabla\times{\bf u})^2{\bf u}+c_7{\bf u}^2{\bf u}\nonumber\\
&+c_8\nabla p+c_9{\bf f}.
\end{align}

The evolution equation for the pressure can be constructed in a similar manner, with $N_p$ that should be a scalar.
The pressure should be a function of the velocity only, so keeping the leading order (in $\nabla$ and ${\bf u}$) term, we will find
\begin{align}\label{eq:p}
\partial_tp &= -\kappa\nabla\cdot{\bf u},
\end{align}
where $\kappa$ is another unknown parameter. 
Using the scaling freedom in defining the latent field $p$ explicit in the equations \eqref{eq:u} and \eqref{eq:p}, without the loss of generality we can set $|c_8|=1$. 
Similarly, we can set $c_9=1$, which amounts to choosing a particular scale for the (unknown) forcing.
The remaining constant (to preserve the translational symmetry in space and time) parameters $c_1$ through $c_7$ (and possibly $\kappa$) need to be determined from data using symbolic regression.

The combination of symmetry and basic physics we used constrains the form of the evolution equations rather significantly, yielding a set of evolution equations with rather few superfluous terms.
Indeed, under certain assumptions, a two-dimensional model of the form
\begin{align}\label{eq:2dNSE}
\partial_t{\bf u} &= -\beta({\bf u}\cdot\nabla){\bf u}+\nu\nabla^2{\bf u}-\alpha{\bf u}
+\rho^{-1}({\bf f}-\nabla p),\\
\label{eq:pc}
\partial_tp &= -\kappa\nabla\cdot{\bf u},
\end{align}
which is a special case of the general model \eqref{eq:u}-\eqref{eq:p}, can be derived analytically for this flow by depth-averaging the three-dimensional Navier-Stokes equation \cite{suri_2014,pallantla_2018}.
Here $\alpha$, $\nu$, $\beta$, $\rho$, and $\kappa$ are constants representing, respectively, the vertical momentum transport, the horizontal momentum transport, the attenuation of inertia due to vertical velocity stratification, the density of the fluid, and the scale of the hydrostatic pressure.

Comparison of numerical simulations and experimental observations over a range of Reynolds numbers $Re$ suggests that the model \eqref{eq:2dNSE}-\eqref{eq:pc} with constant parameters $\alpha$, $\nu$, and $\beta$ is qualitatively accurate \cite{tithof_2017}, but there is a systematic discrepancy that can be mostly accounted for by making these parameters weakly dependent on the Reynolds number $Re$. 
The assumptions made in deriving the model mainly affect the vertical momentum transport represented by the Rayleigh friction term $-\alpha{\bf u}$.
This momentum transport increases with $Re$ due to the advection, which can be accounted for by making the scalar coefficient $\alpha$ velocity-dependent.
This is what the last few terms in the general model \eqref{eq:u} represent. Specifically,
\begin{align}\label{eq:alpha}
\alpha&=-c_3-c_4\nabla\cdot{\bf u}-c_5(\nabla\cdot{\bf u})^2\nonumber\\
&-c_6(\nabla\times{\bf u})^2- c_7{\bf u}^2
\end{align}
can be thought of as a second-order (in ${\bf u}$ and $\nabla$) model of the Rayleigh friction coefficient.

\section{Symbolic regression}
\label{sec:sr}

The data characterizing both components of the velocity field ${\bf u}$ can be obtained, for instance, using particle image velocimetry \cite{suri_2014,tithof_2017} and is assumed to be on a uniform grid $(i,j,k)$, where $i$, $j$, and $k$ correspond to the $x$, $y$, and $t$ directions, respectively.
However, unlike standard symbolic regression problems where all the variables are directly observable, in our problem neither the pressure $p$ nor the forcing ${\bf f}$ are, so both fields have to be either determined independently or eliminated.
In principle, for a fully resolved incompressible ($\kappa\to\infty$) flow field, if ${\bf f}$ were known, $p$ could be obtained in a standard way by applying a divergence to \eqref{eq:u}, which yields a pressure Poisson equation.
Typical experimental data however have a resolution that is too poor (and noise level that is too high) to make it possible to compute pressure in this manner.

In the following, we will focus just on the evolution equation for the velocity field; the evolution equation for the pressure is very simple and the coefficient $\kappa$ can be eliminated altogether by rescaling $c_8$ and/or $p$.
The terms involving both latent fields can be eliminated from \eqref{eq:u} by applying an operator $\hat{P} = \hat{S}\hat{C}\hat{T}$ composed of three operations: $\hat{C}=\hat{z}\cdot\nabla\times$ removes the dependence on $\nabla p$ which is curl-free, $\hat{T}=\partial_t$ removes the dependence on ${\bf f}$ which is constant, and the sparsification operator $\hat{S}$ subsamples the original data in a random fashion.
The corresponding discretization of the resulting PDE (which is second order in time, third order in space, and fourth order overall) has the form  
\begin{equation}\label{eq:symreg}
{\bf q}_0  = Q{\bf c},
\end{equation}
where $Q=[{\bf q}_1\ \dots\ {\bf q}_7]$, ${\bf c}=[c_1\ \cdots\ c_7]^T$ is a vector composed of scalar coefficients to be determined, and the columns
\begin{align}
\begin{matrix*}[l]
{\bf q}_0=\hat{P}(\partial_t{\bf u}-{\bf f}+\nabla p),\quad&
{\bf q}_1=\hat{P}({\bf u}\cdot\nabla){\bf u},\\
{\bf q}_2=\hat{P}\nabla^2{\bf u},&
{\bf q}_3=\hat{P}{\bf u},\\
{\bf q}_4=\hat{P}(\nabla\cdot{\bf u}){\bf u},&
{\bf q}_5=\hat{P}(\nabla\cdot{\bf u})^2{\bf u},\\
{\bf q}_6=\hat{P}(\nabla\times{\bf u})^2{\bf u},&
{\bf q}_7=\hat{P}{\bf u}^2{\bf u} 
\end{matrix*}
\end{align}
correspond to different terms in \eqref{eq:u}. 
Note that ${\bf q}_0=\hat{P}\partial_t{\bf u}$, so none of the terms ${\bf q}_i$ in fact depend on either $p$ or ${\bf f}$.
For the number $K$ of points in the sample exceeding the number of unknown coefficients, this yields an overdetermined system \eqref{eq:symreg} of linear equations for ${\bf c}$, where the ``library'' $Q$ and the ``target'' ${\bf q}_0$ can be evaluated using any algorithm sufficiently robust with respect to noise and sparsity of the data. The particular procedure used in the present paper is described in the next section.

We performed symbolic regression using the iterative algorithm for sparse identification of nonlinear dynamical systems (SINDy) \cite{brunton2016}, which involves computing the solution ${\bf c}$ that minimizes the residual 
\begin{align}\label{eq:eta}
\eta = \| {\bf q}_0 - Q{\bf c} \|_1
\end{align}
of the linear system \eqref{eq:symreg}, followed by a thresholding procedure to remove dynamically irrelevant terms.
Note that the library terms ${\bf q}_i$ can differ by many orders of magnitude. Since it is the product, $c_i{\bf q}_i$, that determines a given term's role in the model, we employ a slightly modified thresholding procedure.
We compare the norms of the products $c_i{\bf q}_i$ to the residual $\eta$: the columns of $Q$ for which $\|c_i{\bf q}_i\|_1<\gamma\eta$ are removed,
and the process is repeated until all remaining terms are above the threshold.
Here $\gamma$ is a constant that can be above or below, but is close to, unity.

\section{Polynomial Interpolation}
\label{sec:pi}

All of the library terms involve spatial and/or temporal derivatives of the velocity field.
Using total variation regularization of the data \cite{rudin1992} to reduce the influence of noise is both prohibitively expensive in higher dimensions and unnecessary given the sparse nature of the system \eqref{eq:symreg}.
Therefore, to accomplish the task of smoothing noisy data and taking numerical derivatives concurrently, a higher-dimensional generalization of the polynomial interpolation in \cite{savitzky1964smoothing} was used instead. 
(We also investigated computation of derivatives using discrete Fourier transform, and found the results to be comparable).
At each point chosen by the sparsification operator $\hat{S}$, the velocity fields were approximated by a polynomial in $x$, $y$, and $t$ fitted to discrete data on a rectangular domain $\Omega$ of size $2H_x\times 2H_y\times 2H_t$ centered at a grid point $(x_i,y_j,t_k)$. 

In particular, the $x$-component of the velocity $u(x,y,t)$ near $(x_i,y_j,t_k)$ was approximated as
\begin{equation}
\tilde{u}(x,y,t) =  \sum^L_{l=0}\sum^M_{m=0}\sum^N_{n=0}U_{ijk}^{lmn}\bar{x}^l\bar{y}^m\bar{t}^n,
\label{eq:u_poly}
\end{equation}
where the overbar denotes the shifted and rescaled coordinates in which the domain $\Omega$ becomes a cube $\Omega'=[-1,1]\times[-1,1]\times[-1,1]$, e.g.,
\begin{equation}
\bar{x} = \frac{x-x_i}{H_x},
\end{equation}
etc. 
The order of the polynomial in each direction should be at least as large as the order of the highest derivative appearing in the model equation \eqref{eq:u} after the operator $\hat{C}\hat{T}$ is applied, but ultimately is a tunable parameter, with the specific choice to be discussed in more detail in the next section.

The coefficients $U_{ijk}^{lmn}$ were found by minimizing the cost function
\begin{equation}
F = \sum_{(i,j,k)\in\Omega} w_{ijk} (u_{ijk} - \tilde{u}(x_i,y_j,t_k))^2,
\label{eq:F}
\end{equation}
where $w_{ijk}$ is a weighting function.
This is a standard least squares problem whose solution is given by setting $\partial F/\partial U_{ijk}^{qrs} = 0$. This yields a system of $(L+1)(M+1)(N+1)$ linear equations
\begin{equation} \label{eq:linpoly}
\langle w\:u\:\bar{x}^{q}\bar{y}^{r}\bar{t}^{s}\rangle_{\Omega} = \sum^{L,M,N}_{l,m,n}U_{ijk}^{lmn} \langle w\:\bar{x}^{l+q}\bar{y}^{m+r}\bar{t}^{n+s}\rangle_{\Omega},
\end{equation}
where $\langle\cdot\rangle_{\Omega}$ denotes the average over the spatiotemporal sub-domains for which the local fits are defined. The weighting function $w_{ijk}$ was used to bias the accuracy of the approximation toward the central point of the domain $\Omega$ (where all of the derivatives are evaluated) and is defined as a Gaussian
\begin{equation}
w_{ijk} = \exp\left(-\frac{\bar{x}^2+\bar{y}^2+\bar{t}^2}{\lambda^2}\right),
\end{equation}
with the width $\lambda$ being another tunable parameter of the model (we set $\lambda=0.5$). The same procedure was used to determine the coefficients $V_{ijk}^{lmn}$ for the $y$-component of the velocity $v(x,y,t)$.

After the polynomial coefficients have been determined, the row of the library $Q$ and the target ${\bf q}_0$ corresponding to the point $(x_i,y_j,t_k)$ can be constructed by evaluating the respective derivatives of $u$ and $v$ at $(\bar{x},\bar{y},\bar{t})=(0,0,0)$ using \eqref{eq:u_poly}.
For instance, 
\begin{align}
q_2^{ijk} = 6V_{ijk}^{301} + 2V_{ijk}^{121} - 2U_{ijk}^{211} - 6U_{ijk}^{031}.
\end{align}
The process was repeated for each point defined by $\hat{S}$ in order to completely evaluate the library and the target.
$\hat{S}$ was defined by randomly selecting the points on the entire 3D grid representing the spatially and temporally discretized trajectory.
Throughout the paper, $K=250$ points were used to construct the library; 
neither the mean nor the standard deviation of the coefficients $c_i$ were found to exhibit meaningful variation for a larger number of points.

\section{Results}
\label{sec:r}

Surrogate data used for testing the symbolic regression procedure was generated using the model \eqref{eq:u}-\eqref{eq:p} with the parameters $c_1=-0.826$, $c_2=0.0487$, $c_3=-0.157$, $c_4=0.164$, $c_5=c_6=c_7=0$, $c_8=-1$, and
$\kappa=2015$. 
This set of parameters describes a nearly incompressible flow found in the experiment described in Ref. \cite{tithof_2017}, which features a forcing field with a sinusoidal profile in the $y$ direction with period $2\chi=2$ and amplitude equal to $1.0649$ in nondimensional units.
The solution describing a weakly turbulent flow was obtained using a numerical integration scheme based on operator splitting as described in Ref. \cite{pallantla_2018}. 
The linear terms were evolved in time implicitly, while the nonlinear terms were handled via a $2^{nd}$ order Adams-Bashforth scheme.
The solution was integrated on a computational grid with $\Delta x_c=\Delta y_c=0.025$, and $\Delta t_c\approx0.02$.
Gaussian random noise with variance $\sigma$ was added to both components of the flow velocity ${\bf u}$. For reference, the maximal flow velocity is $O(1)$ in nondimensional units.

In order for the algorithm to produce meaningful results, its various tunable parameters must be properly set. 
The noiseless case exhibits the least amount of sensitivity to variation of fitting parameters; the only restriction is that the polynomial orders $L$, $M$, and $N$ be high enough to capture the variation in the data over $\Omega$. While higher order interpolation allows better approximation of the data, it is also more sensitive to noise.  To mitigate the influence of noise, a larger number of measurements can be used. There are two ways to achieve this: by increasing the size of the sampling domain  $\Omega$ or by using a finer grid on which data are measured.  The largest size of $\Omega$ is effectively limited by the characteristic length and time scales for the problem. 
In the present problem, the natural length scale is defined by $\chi$.
Consequently, we will set $H_x=H_y=\chi/2$, such that the width of $\Omega$ in both spatial dimensions is equal to $\chi=1$.
There is no natural time scale, so we will choose one based on the autocorrelation time $\tau\approx9.9$.  
In the following we set $H_t\approx0.85\tau$ which is an optimal choice for $\sigma=10^{-3}$
and $M=L$.
Furthermore we use the finest grid available to us in space, i.e., $\Delta x=\Delta x_c$, while in time we use $\Delta t=25\Delta t_c$. 
With this choice, $\Omega$ corresponds to a $40\times 40\times 34$ block of data with dimensions that are roughly comparable in the spatial and temporal directions. 
At higher grid sizes evaluating the averages in (\ref{eq:linpoly}) becomes computationally expensive.

\begin{figure}[t!]
\centering
\includegraphics[width=1\linewidth]{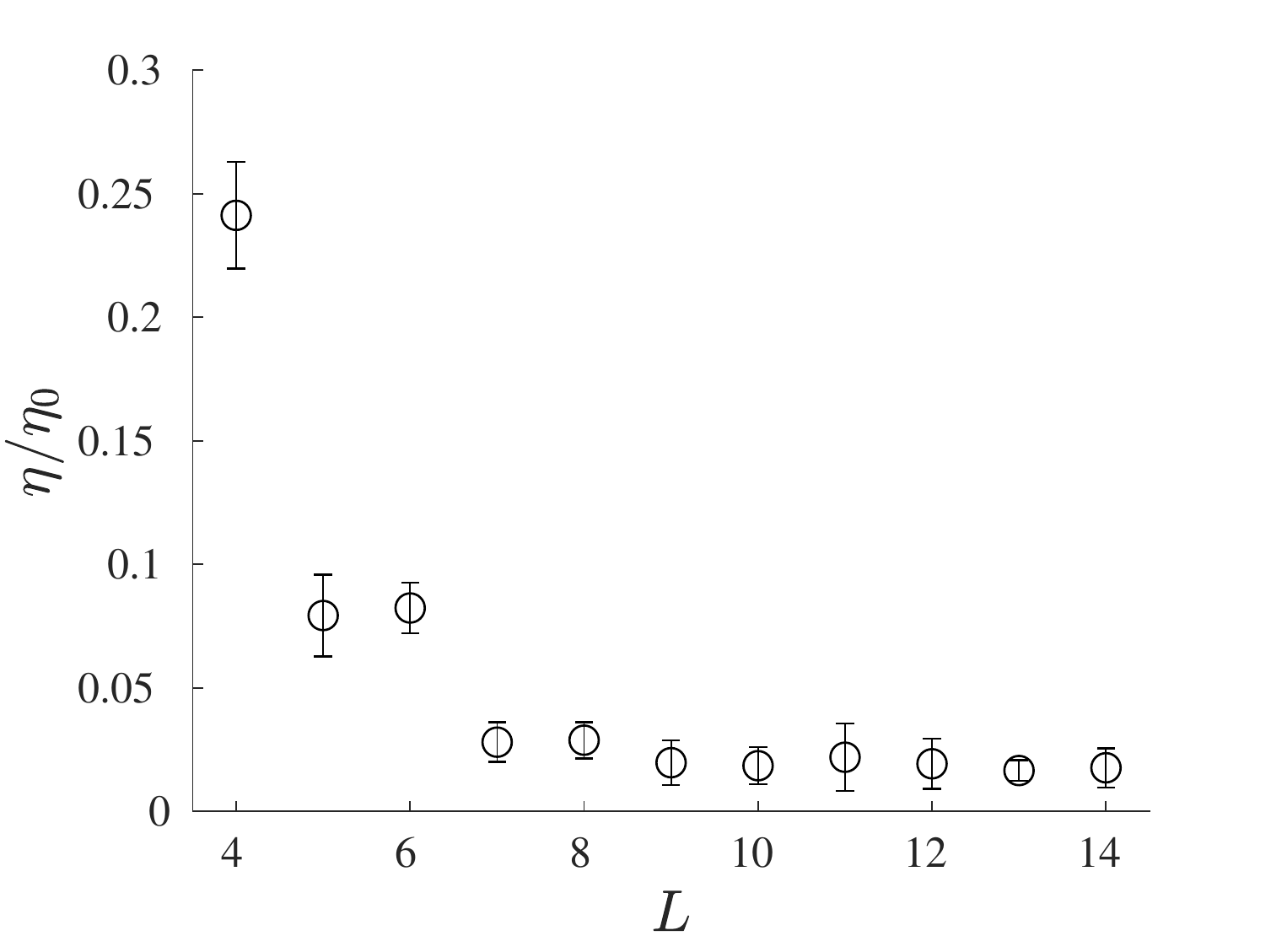}
\caption{Residual as a function of polynomial order for 
$N=10$ and $\sigma=0$. Here, error bars denote standard deviation and symbols denote mean values.}
\label{fig:L_vs_eta_n0}
\end{figure}

To investigate how the choice of polynomial order affects the accuracy of the fit and hence the accuracy with which various partial derivatives of ${\bf u}$ are evaluated, we computed the residual \eqref{eq:eta}.
The dependence of $\eta$ (normalized by the magnitude of the target $\eta_0=\|{\bf q}_0\|_1$) is shown in Fig. \ref{fig:L_vs_eta_n0}. 
Here and below, the averages and standard deviation are computed using an ensemble of $40$ different realizations of the sampling operator $\hat{S}$.
Note that $\eta$
 generally does not vanish even for the noiseless perfect model of the problem due to discretization errors of the numerical solution. 
Also note that the magnitude of $\eta$ describes the accuracy with which equation \eqref{eq:symreg} is satisfied, not the accuracy of the numerical solution to the model \eqref{eq:u}-\eqref{eq:p}.
As expected, $\eta$ decreases for low $L$, but beyond some threshold (in this case $L=7$), increasing the polynomial order has little effect on the residual. In particular, $L=N=10$ results in both a low value of the residual and a small error in parameter estimation in the noiseless case, as we will see below.

\begin{figure}[t!]
\centering
\includegraphics[width=1\linewidth]{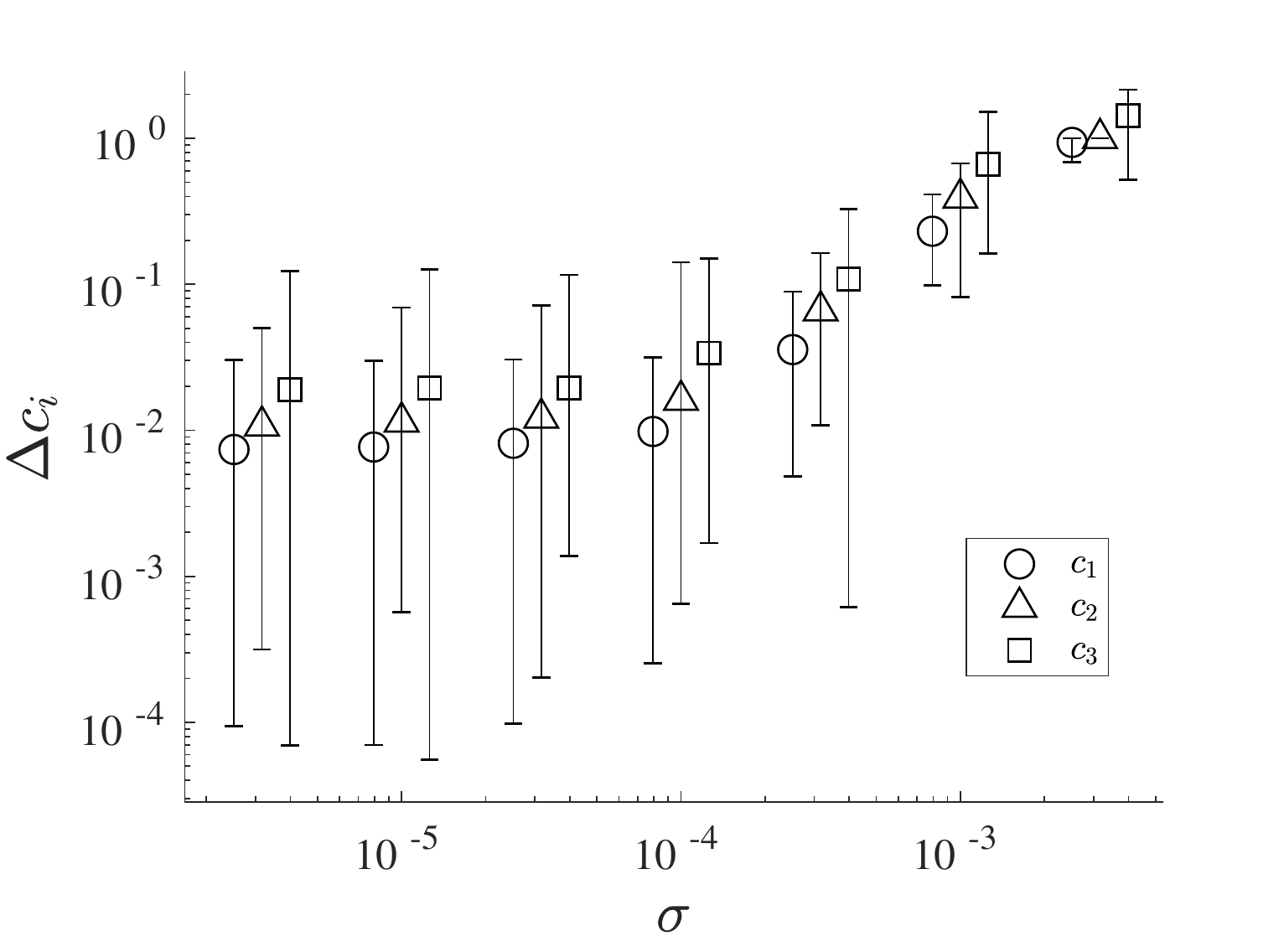}
\caption{Parameter error as a function of noise amplitude for $L=N=10$. The error for $c_4$ is not shown because SINDy discards the corresponding term. Here and below markers are shifted left or right to avoid overlap. Error bars indicate the full range of data, and markers indicate mean values.}
\label{fig:sig_vs_par_n0}
\end{figure}

To determine how the results depend on the amplitude $\sigma$ of measurement noise, we performed symbolic regression 
and compared the coefficients $\tilde{c}_i$ produced by SINDy with the reference values $c_i$ used to generate the surrogate data for $\kappa=2015$. 
Of the four nonzero parameters used in generating the data, three ($c_1$, $c_2$, and $c_3$) were correctly identified as being nonzero and estimated with a small relative error
\begin{align}
\Delta c_i = \left| \frac{c_i - \tilde{c}_i}{c_i} \right|
\end{align}
(of order one percent) for sufficiently small 
$\sigma$, as illustrated by Fig. \ref{fig:sig_vs_par_n0}.
However, the coefficient $c_4$ was incorrectly set to zero by the algorithm for all $\sigma$.
Furthermore, the accuracy in estimating all of the remaining parameters decreased sharply for $\sigma\gtrsim 10^{-4}$.

The failure of symbolic regression to correctly identify the value of $c_4$ can be understood qualitatively by recalling that it is the product $c_4{\bf q}_4$ whose magnitude is used to determine whether the corresponding term should be retained or discarded.
For our choice of parameters, $\|c_4{\bf q}_4\|_1\lesssim\eta$, 
suggesting that the magnitude of this term is as small or smaller than the accuracy to which the governing equation can be satisfied.
As mentioned previously, the parameter set used here corresponds to a nearly incompressible flow where $\nabla\cdot{\bf u}$ is nonzero but very small.
Indeed, from \eqref{eq:u} and \eqref{eq:p} we find ${\bf q}_4\sim\nabla\cdot{\bf u}\sim \kappa^{-1}\approx 5\times 10^{-4}$.
By eliminating the term $c_4(\nabla\cdot{\bf u}){\bf u}$ representing the effect of compressibility \cite{pallantla_2018} from the model \eqref{eq:u}, symbolic regression effectively recognized this fact.

To verify that this term was indeed eliminated due to the large value of $\kappa$ (and not some shortcoming of the method), we repeated the analysis, setting $\kappa=1$ to amplify the compressibility effects. In this case the term $c_4(\nabla\cdot{\bf u}){\bf u}$ is retained in the model and the value of the parameter $c_4$ is determined correctly 
for sufficiently small 
$\sigma$.
This is a good example illustrating when symbolic regression fails to identify terms that are generally required by the physics of the problem but may be neglected under some conditions.

To make the argument more quantitative, let us introduce the measure 
\begin{align}
R_i=\frac{\|c_i{\bf q}_i\|_1}{\eta},
\end{align}
of the magnitude of a particular term in the linear equation (\ref{eq:symreg}) relative to the corresponding residual \eqref{eq:eta}.
Symbolic regression can correctly identify a particular term in the model only if the corresponding $R_i>1$; furthermore, we can expect the accuracy of parameter reconstruction to decrease as $R_i$ approaches unity. For our choice of fitting parameters, $R_4$ is below unity for $\kappa=2015$ and above unity ($R_4\approx20$) for $\kappa=1$.
Correspondingly, the terms with lower $R_i$ exhibit the worst fitting accuracy; this explains the larger relative error in $c_3$ compared with $c_1$ and $c_2$ in the $\kappa=2015$ case (cf. Fig. \ref{fig:sig_vs_par_n0}), since although $R_3\approx30$, it is much smaller than $R_1$ and $R_2$, which both have $R_i>100$.

\begin{figure}[t!]
\centering
\includegraphics[width=1\linewidth]{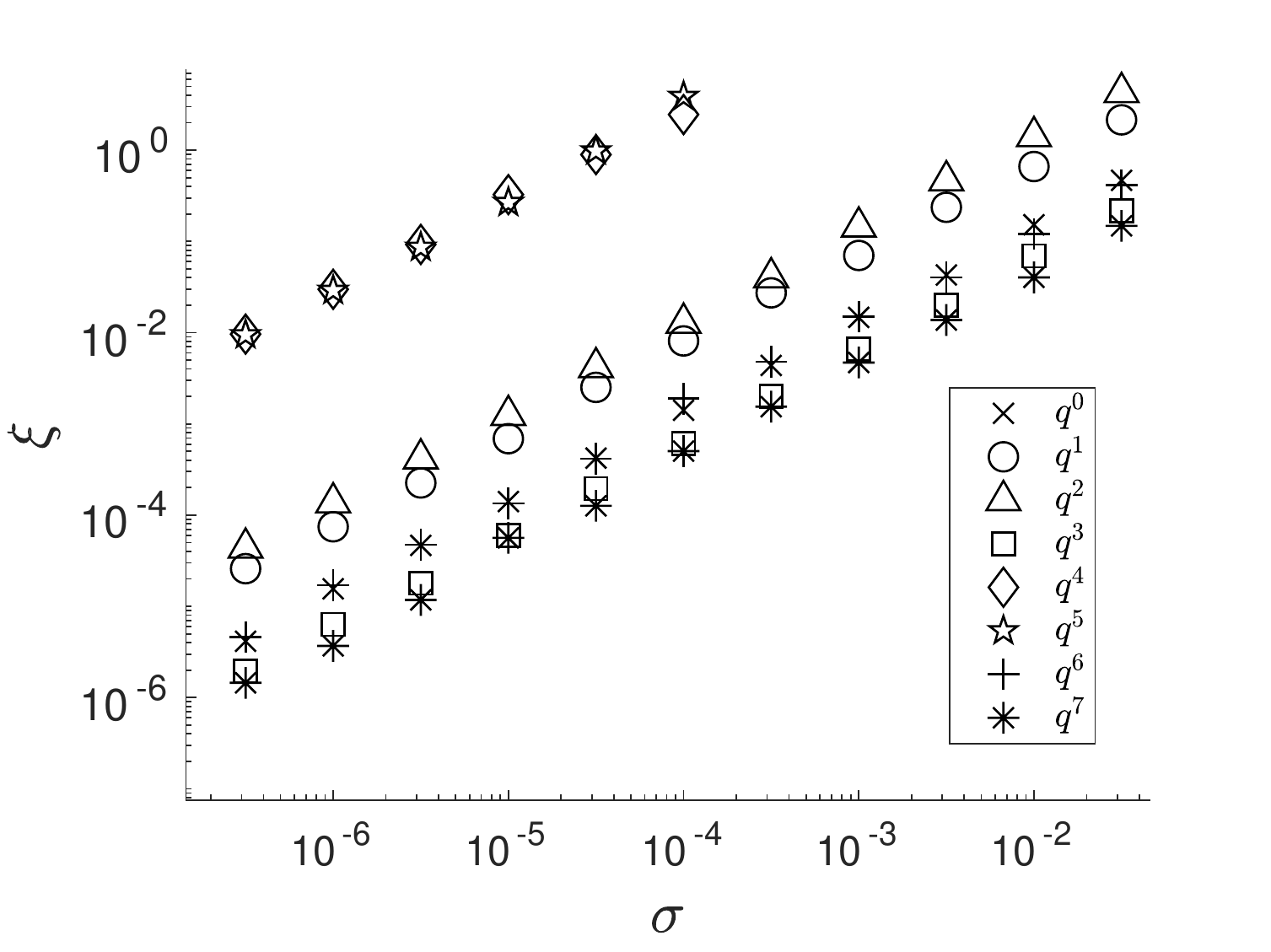}
\caption{Relative accuracy of different library terms as a function of noise amplitude for  $L=N=10$.
}
\label{fig:sig_vs_xi_n0}
\end{figure}

Figure \ref{fig:sig_vs_par_n0} also shows that, for the choice of fitting parameters optimized for noiseless data, the accuracy of symbolic regression sharply decreases for $\sigma>O(10^{-4})$. To understand why this happens, let us define the relative accuracy with which a particular library term is evaluated over the entire sample
\begin{align}
\xi_i(\sigma)= \frac{\|{\bf q}_i(0) - {\bf q}_i(\sigma)\|_\infty }{ \|{\bf q}_i(0)\|_\infty}.
\end{align}
The effect of noise on the accuracy of all the library terms is shown in Fig. \ref{fig:sig_vs_xi_n0}. Note that the lowest accuracy (highest $\xi_i$) corresponds to the terms ${\bf q}_4$ and ${\bf q}_5$ which are linear and quadratic, respectively, in $\nabla\cdot\bf u$, which is very small. These terms are the most susceptible to corruption by noise but, for large values of $\kappa$, they are eliminated by symbolic regression anyway. As might be expected, in the absence of these two terms, the term ${\bf q}_2$, which involves the highest order derivative (third order in space and first in time), is the least accurate in the presence of noise. This helps explain the difficulties symbolic regression has with identifying high order derivatives in all PDE models in the presence of noise. For instance, in a previous study \cite{rudy2017}, the coefficient of the fourth order derivative term in the Kuramoto-Sivashinsky was determined with a 52\% error in the presence of 1\% noise. In our case, the terms ${\bf q}_3$ and ${\bf q}_7$ which involve the lowest order derivative (first in space and time), have the smallest error, suggesting that the order of the derivative is one of the main factors which determine the accuracy of regression in the presence of noise.

\begin{figure}[t!]
\centering
\includegraphics[width=1\linewidth]{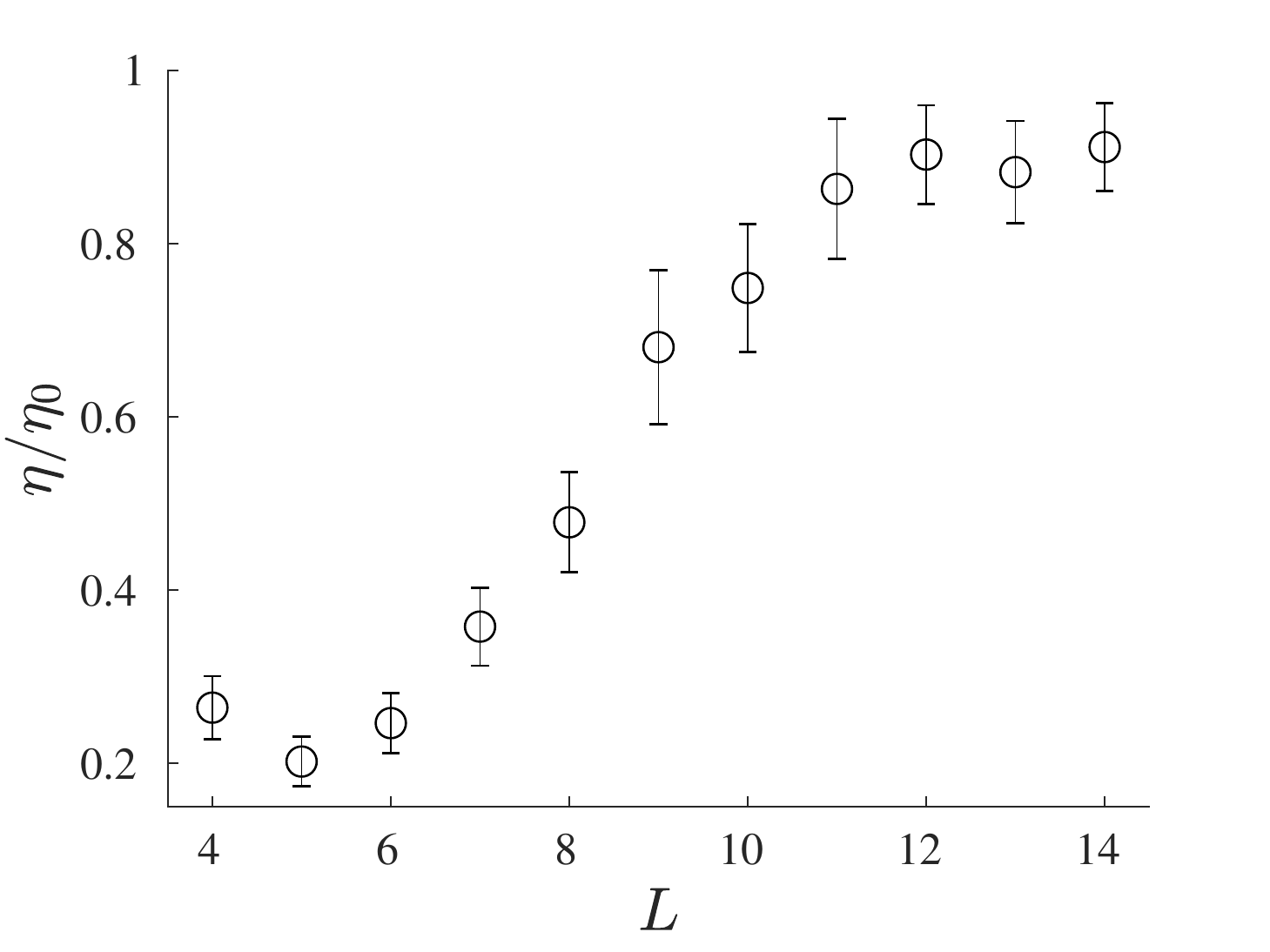}
\caption{Residual as a function of polynomial order for 
$N=10$ and $\sigma=10^{-3}$. Here, error bars denote standard deviation and symbols denote mean values.}
\label{fig:L_vs_eta_n3}
\end{figure}

The effect of noise can be offset, to some extent by a different choice of parameters. In particular, the order of the polynomial interpolation can be reduced to decrease noise sensitivity. The dependence of the residual $\eta$ on $L$ is shown in Fig. \ref{fig:L_vs_eta_n3} for noise amplitude $\sigma=10^{-3}$ at which our previous choice of parameters lead to unacceptably large errors. At this value of $\sigma$, we find a minimum around $L=5$ (with a significant increase in $\eta$ compared to the noiseless case), which represents a balance between the accuracy of the interpolation in capturing the spatial variation of the data at higher $L$ and the noise sensitivity at lower $L$. In fact, we found that setting $L=M=N=6$ is the best choice for minimizing both the residual and the error in parameter estimation.

\begin{figure}[t!]
\centering
\includegraphics[width=1\linewidth]{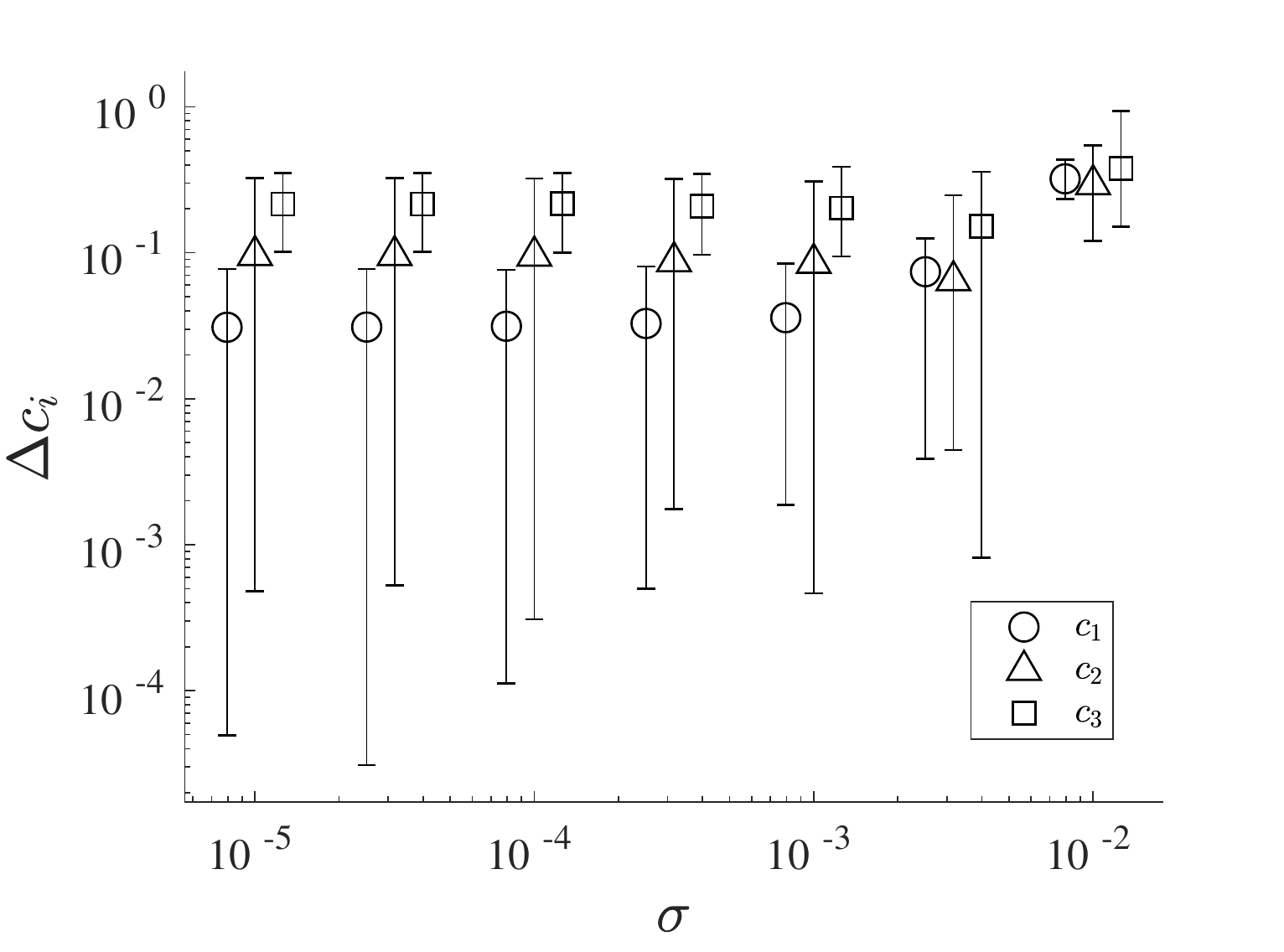}
\caption{Parameter error as a function of the noise amplitude for $L=N=6$. The error for $c_4$ is not shown because SINDy discards the corresponding term. Error bars indicate the full range of data, and markers indicate mean values.}
\label{fig:sig_vs_par_n3}
\end{figure}

Using the fitting parameters optimized for higher noise levels, symbolic regression identifies the correct model (aside from the negligible term ${\bf q}_4$) 
with all the model parameters estimated to within $\sim$10\% for 0.3\% noise and to within $\sim$30\% for 1\% noise, as illustrated by Fig. \ref{fig:sig_vs_par_n3}. This is comparable to the accuracy achieved for the Kuramoto-Sivashinsky equation \cite{rudy2017}, which also includes a fourth order derivative. The trade-off of this choice of fitting parameters is the decrease in the accuracy of model parameter estimation at lower noise levels. Furthermore, $R_3$ becomes close to unity, so symbolic regression yields false negatives for a noticeable fraction of the trials (the data shown in  Fig. \ref{fig:sig_vs_par_n3} was calculated after discarding the results for which SINDy eliminated the term ${\bf q}_3$). In comparison, false negatives did not appear for $L=N=10$ until fairly high levels of noise. 
These false negatives occur for lower $L$ and/or $N$ because
the magnitude of the residual is determined by the error in the term(s) most affected by the insufficiently accurate approximation (here, the term ${\bf q}_2$ which involves the highest order derivative).
For the lower $L$, the variation in the data is not fully resolved, meaning that $R_3$ is pushed closer to unity (and hence $c_3$ can be estimated with less accuracy).

\section{Conclusions}
\label{sec:c}

This paper presents an approach that allows symbolic regression to be extended for data-driven discovery of PDE-based models which involve unobservable and/or unknown (latent) variables. Our approach relies on two key ideas: (1) using spatial symmetry as well as other physical constraints to select the terms that can appear in the model and (2) applying a differential operator designed to remove terms which involve the latent variables. We illustrated these ideas by using sparse regression to construct a two-dimensional model for a weakly compressible Kolmogorov-like flow in a thin electrolyte layer driven by a steady Lorentz force with an unknown spatial profile. In this particular case, two latent variables -- the forcing field which is a vector and the pressure which is a scalar -- have been eliminated by applying, respectively, a temporal derivative and a curl to a nonlinear model with nine different terms allowed by symmetry. 

While previous studies have demonstrated the power of sparse regression for data-driven model discovery, they left a number of questions unanswered. In particular, how should one choose the threshold that determines which terms in the model are relevant? Since those studies mainly focused on reconstructing well-known models, the threshold could be chosen in an ad hoc fashion such that the {\it a priori} known model was recovered. In case the form of the model is not known {\it a priori}, the proper choice is less clear, since sparse regression will recover different models for different choices of the threshold. We have shown that a self-consistent choice should be based on the residual $\eta$ of the linear system \eqref{eq:symreg}: in most cases symbolic regression can be considered to have successfully reconstructed the model once the corresponding norm of every remaining term $c_i{\bf q}_i$ is larger than the residual. 

However, while the residual is the proper metric for determining the relevance of different terms, fine-tuning the threshold 
does have an effect on the reconstructed model, mainly affecting the terms with $R_i=O(1)$. 
For $\gamma$ somewhat smaller than unity, symbolic regression can produce false positives, e.g.,
the term $c_4(\nabla\cdot{\bf u}){\bf u}$ is retained even in the essentially incompressible case when $\nabla\cdot{\bf u}$ is very small and $R_4$ is just below unity. However, the value (and sign) of the coefficient $c_4$ is found to vary drastically for different realizations of the sampling operator $\hat{S}$, suggesting that this term is not dynamically relevant.
For $\gamma$ somewhat larger than unity, this term is correctly removed by the thresholding procedure for sufficiently low noise. 
However, this choice also causes the term $c_3{\bf u}$ to be removed for some realizations of $\hat{S}$ at higher noise amplitudes, when $R_3$ becomes comparable to unity as well.
In case a relevant term (such as $c_3{\bf u}$) is removed by the thresholding procedure, the residual increases noticeably, which allows detection of false negatives.
Manually including the term $c_3{\bf u}$ in such cases decreases the residual by about 10\%, indicating that it is dynamically relevant despite being smaller than the residual.
Fine-tuning $\gamma$ based on these metrics can allow more robust results in border-line cases.

Our study has also highlighted the major weakness of all spatially local approaches to sparse regression. Regardless of whether one uses a polynomial interpolation of the data, total variation regularization, or some other similar approach to construct the linear system whose solution determines the model, the procedure is inherently sensitive to noise, especially when higher order derivatives are involved. This difficulty is well-illustrated by both the failure of previous studies to reconstruct with acceptable accuracy the fourth-order Kuramoto-Sivashinsky equation in the presence of as little as 1\% noise and a similar loss in the accuracy for the model considered here, whose latent-variable-free form also involves fourth order derivatives. In both cases, as the noise amplitude increases, the term involving the highest order derivative becomes the largest contributor to the residual at which point its coefficient cannot be reliably determined anymore. Since experiments commonly involve substantially higher amounts of noise, a more robust alternative to such local methods is needed.

\begin{acknowledgments}

This material is based upon work supported by the National Science Foundation under Grant No. CMMI-1725587.

\end{acknowledgments}

\bibliographystyle{unsrt}
\bibliography{./latent_discovery}

\end{document}